\documentclass[12pt]{scrartcl}
\usepackage[american]{babel}
\usepackage{amsthm,amsmath,amssymb}
\newcommand*{\set}[1]{\{#1\}}
\newcommand*{\defset}[2]{\{#1\,\mid\,#2\}}
\newcommand{\OO}{\mathcal O}
\theoremstyle{plain}
\newtheorem{theorem}{Theorem}
\usepackage{stix2}
\usepackage{xcolor}
\usepackage[colorlinks=true,linkcolor=blue,citecolor=purple]{hyperref}
\begin{document}
\title{Another proof of Segre's Theorem about ovals}
\author{Peter M\"uller}
\maketitle
\abstract{\noindent\textbf{Abstract:} In 1955 B.~Segre showed that any
  oval in a projective plane over a finite field of odd order is a
  conic. His proof constructs a conic which matches the oval in some
  points and tangents, and then shows that it actually coincides with
  the oval. The different proof given here parametrizes an affine
  piece of the oval and shows directly that the parametrization is
  given by a polynomial of degree $2$.
  \par\vskip\baselineskip\noindent
\par\vskip\baselineskip\noindent
}
\section{Introduction} An oval in a finite projective plane of order
$q$ is a set of $q+1$ points such that no three of them are on a line.

It is easy to see that each non-singular conic in the projective plane
$\operatorname{PG}(2,\mathbb F_q)$ is an oval. By Segre's well known
theorem \cite{Segre:Ovals} there are no other ovals if $q$ is odd. His
proof consists of two parts: If $A$ is a point on the oval $\OO$, then
by the definition of an oval there is a unique line $t_A$ through $A$,
called the tangent in $A$, which intersect $\OO$ only in $A$. Let
$A_0,A_1,A_2$ be distinct points on $\OO$. In the first step he shows
that the $3$-tangents degeneration of Brianchon's theorem holds for
$\OO$: Let $l_i$ be the line through $A_i$ and
$t_{A_{i-1}}\cap t_{A_{i+1}}$ (indices taken modulo $3$), then
$l_1,l_2,l_3$ intersect in a single point. This is used in the second
step to show that the conic which shares three points of $\OO$ and
their tangents in these points actually coincides with $\OO$.

His proof went essentially unchanged into many sources, see e.g.\
\cite[Section 9.7]{Cameron:Combinatorics}, \cite[Theorem
8.14]{Hirschfeld:PG}, and \cite[6.~Appendix]{HughesPiper}. A more
expository presentation is in \cite{Browne_etal:Segre}. The lecture
notes \cite{Hartmann:PlanarCircleGeometries} contain a variation of
Segre's argument. The first step verifies the $3$-point degeneration
of Pascal's theorem: Let $P_i$ be the intersection of $t_{A_i}$ with
the line through $A_{i-1}$ and $A_{i+1}$. Then the points
$P_0, P_1, P_2$ are collinear. The second step then shows that this
again implies that $\OO$ is a conic.

The proof of Segre's theorem in this note uses ideas similar to
Segre's, while it differs in some details. Instead of constructing a
conic and showing that it coincides with the given oval, we
parametrize an affine piece of the oval and show that the
parametrization is given by a quadratic polynomial.
\begin{theorem}\label{T}
  Let $F$ be a finite field of odd order, and $f:F\to F$ be a map such
  that no three points of the graph $\defset{(x,f(x))}{x\in F}$ are
  collinear in the affine plane over $F$. Then $f$ is given by a
  quadratic polynomial.
\end{theorem}
Theorem \ref{T} had been used previously in a non-geometric context,
where it was derived as an easy consequence of Segre's theorem. See
e.g.\ \cite[p.~99]{Gluck:PlanarPol} and \cite[Result 3]{RonyaiSzonyi}.

The purpose of this note is to go the other way around. We give a
direct proof of Theorem \ref{T}, and show how Segre's theorem follows
from that.
\section{Theorem \ref{T} implies Segre's theorem}
Let $\OO$ be an oval in the projective $(X:Y:Z)$--plane over the
finite field $F$ of odd order. After a suitable collineation we may
and do assume that $(0:1:0)\in\OO$, and that $Z=0$ is the tangent line
through this point. Since each line $X=xZ$ ($x\in F$) contains
$(0:1:0)$ and is not the tangent, we obtain the following: For each
$x\in F$ there is a unique $y\in F$ such that $(x:y:1)$ is on the
oval. Thus the hypothesis of Theorem \ref{T} holds for the map
$f:F\to F,x\mapsto y$, hence $f(x)=ax^2+bx+c$ for $a,b,c\in F$ and all
$x\in F$. Therefore $\OO$ is the conic given by $YZ=aX^2+bXZ+cZ^2$.
\section{Proof of Theorem \ref{T}}
  Set $q=\vert F\rvert$ and fix $u\in F$. The $q-1$ slopes
  $\frac{f(x)-f(u)}{x-u}$, $x\in F\setminus\set{u}$, are pairwise
  distinct elements from $F$. Let $s(u)\in F$ denote the unique
  element which is not among these slopes, so
\[
  \defset{\frac{f(x)-f(u)}{x-u}}{x\in F\setminus\set{u}}=%
  F\setminus\set{s(u)}.
\]
Now fix $v\in F$ different from $u$. Then
\[
\defset{\frac{f(x)-f(u)}{x-u}}{x\in F\setminus\set{u,v}}=%
F\setminus\set{s(u),\frac{f(v)-f(u)}{v-u}}.
\]
Upon replacing each element $t$ in these sets by
$t-\frac{f(v)-f(u)}{v-u}$ we obtain
\begin{equation}\label{sets}
\defset{\frac{g(x,u,v)}{x-u}}{x\in F\setminus\set{u,v}}=%
F\setminus\set{s(u)-\frac{f(v)-f(u)}{v-u}, 0},
\end{equation}
where
\begin{equation}
  g(x,u,v)=f(x)-f(u)-(x-u)\frac{f(v)-f(u)}{v-u}.
\end{equation}
Let $P$ be the product of the nonzero elements of $F$. Taking the
product of the elements in the sets of equation \eqref{sets}, and
noting that $\prod_{x\in F\setminus\set{u,v}}(x-u)=P/(v-u)$, we
obtain
\begin{equation}\label{prod}
\frac{\prod_{x\in F\setminus\set{u,v}}g(x,u,v)}{P/(v-u)}=
\frac{P}{s(u)-(f(v)-f(u))/(v-u)}.
\end{equation}
Note that
\[
  g(x,u,v)-g(x,v,u)=f(v)-f(u)-(x-u-x+v)\frac{f(v)-f(u)}{v-u}=0,
\]
so the left hand side of \eqref{prod} changes sign upon switching $u$
and $v$, and therefore so does the right hand side.
We obtain
\[
2(f(u)-f(v))=(s(u)+s(v))(u-v)
\]
for all $u,v\in F$. Setting $v=1$ and $v=0$ we get
\begin{align*}
  2(f(u)-f(1)) &= (s(u)+s(1))(u-1)\\
  2(f(u)-f(0)) &= (s(u)+s(0))u.
\end{align*}
Subtract $u-1$ times the second equation from $u$ times the first one
to obtain
\[
  2f(u) = (s(1) - s(0))u(u-1) + 2(f(1) - f(0))u + 2f(0).
\]
As $q$ is odd, we see that $f$ is given by a polynomial of degree
$\le2$. The theorem follows, since on the other hand $f$ is neither
linear nor constant.

\noindent{\scshape Institut f\"ur Mathematik, Universit\"at W\"urzburg,
  Campus Hubland Nord, 97074 W\"urzburg, Germany}\par
\noindent{\slshape E-mail: }\texttt{peter.mueller@uni-wuerzburg.de}

\begin{thebibliography}{1}
\providecommand{\url}[1]{{#1}}
\providecommand{\urlprefix}{URL }
\expandafter\ifx\csname urlstyle\endcsname\relax
  \providecommand{\doi}[1]{DOI~\discretionary{}{}{}#1}\else
  \providecommand{\doi}{DOI~\discretionary{}{}{}\begingroup
  \urlstyle{rm}\Url}\fi

\bibitem{Browne_etal:Segre}
Browne, P.J., Dougherty, S.T., \'O{}~Cath\'ain, P.: Segre's theorem on ovals in
  {D}esarguesian projective planes.
\newblock Irish Math. Soc. Bull. \textbf{{}}(91), 37--47 (2023).
\newblock \doi{10.33232/bims.0091.37.47}.
\newblock \urlprefix\url{https://doi.org/10.33232/bims.0091.37.47}

\bibitem{Cameron:Combinatorics}
Cameron, P.J.: Combinatorics: topics, techniques, algorithms.
\newblock Cambridge University Press, Cambridge (1994)

\bibitem{Gluck:PlanarPol}
Gluck, D.: A note on permutation polynomials and finite geometries.
\newblock Discrete Math. \textbf{80}(1), 97--100 (1990).
\newblock \doi{10.1016/0012-365X(90)90299-W}.
\newblock \urlprefix\url{https://doi.org/10.1016/0012-365X(90)90299-W}

\bibitem{Hartmann:PlanarCircleGeometries}
Hartmann, E.: Planar {C}ircle {G}eometries (2013).
\newblock
  \urlprefix\url{https://www2.mathematik.tu-darmstadt.de/~ehartmann/circlegeom.pdf}.
\newblock Lecture Notes

\bibitem{Hirschfeld:PG}
Hirschfeld, J.W.P.: Projective geometries over finite fields, second edn.
\newblock Oxford Mathematical Monographs. The Clarendon Press, Oxford
  University Press, New York (1998)

\bibitem{HughesPiper}
Hughes, D.R., Piper, F.C.: Projective planes, \emph{Graduate Texts in
  Mathematics}, vol. Vol. 6.
\newblock Springer-Verlag, New York-Berlin (1973)

\bibitem{RonyaiSzonyi}
R\'onyai, L., Sz{\H o}nyi, T.: Planar functions over finite fields.
\newblock Combinatorica \textbf{9}(3), 315--320 (1989).
\newblock \doi{10.1007/BF02125898}.
\newblock \urlprefix\url{https://doi.org/10.1007/BF02125898}

\bibitem{Segre:Ovals}
Segre, B.: Ovals in a finite projective plane.
\newblock Canadian J. Math. \textbf{7}, 414--416 (1955).
\newblock \doi{10.4153/CJM-1955-045-x}.
\newblock \urlprefix\url{https://doi.org/10.4153/CJM-1955-045-x}

\end{thebibliography}
\end{document}